\documentclass{article}
\usepackage[utf8]{inputenc}
\usepackage[T1]{fontenc}
\usepackage[margin=1.2in]{geometry}
\usepackage{lmodern}
\usepackage{subcaption}
\usepackage{mathtools}
\usepackage{amsthm,url}
\usepackage{todonotes}
\usepackage{mathrsfs,amssymb} 

\usepackage{enumitem}

\newtheorem{theorem}{Theorem}

\newtheorem{lemma}[theorem]{Lemma}

\newtheorem{conj}[theorem]{Conjecture}

\theoremstyle{definition}
\newtheorem*{definition*}{Definition}
\newtheorem{ex}[theorem]{example}

\usepackage{pgfmath, pgfplots}
\pgfplotsset{width=10cm,compat=1.9}
\newcommand*{\myproofname}{Proof}

\newcommand{\pr}{\operatorname{Pr}}

\def\eps{\varepsilon}
\def\F{\mathcal{F}}
 
 \renewcommand{\Pr}{\,\mathbb{P}}
 \newcommand{\Ss}{\mathcal{P}}
\newcommand*{\abs}[1]{\lvert #1\rvert}

\title{Progress on the union-closed conjecture and offsprings in winter 2022-2023}
\author{Stijn Cambie\thanks{Extremal Combinatorics and Probability Group (ECOPRO), Institute for Basic Science (IBS), Daejeon, South Korea, supported by the Institute for Basic Science (IBS-R029-C4), E-mail: {\tt stijn.cambie@hotmail.com}} }

\begin{document}

\maketitle

Mathematicians had little idea whether the easy-to-state union-closed conjecture was true or false even after $40$ years. However, last winter saw a surge of interest in the conjecture and its variants, initiated by the contribution of a researcher at Google. Justin Gilmer made a significant breakthrough by discovering a first constant lower bound for the proportion of the most common element in a union-closed family.

\section{Introduction of the Union-Closed conjecture}

The union-closed conjecture is due to Peter Frankl\footnote{See also \url{https://en.wikipedia.org/wiki/P\%C3\%A9ter_Frankl} and \url{https://www.nrc.nl/nieuws/2023/01/20/na-wiskundige-opwinding-op-sociale-media-is-het-probleem-van-de-kleurige-knikkers-bijna-opgelost-2-a4154586}.}, who constructed the elegant statement in $1979$ after observing many implications of the statement.
Before fully stating it, we need to define crucial concepts from set theory.

The ground set is generally denoted with $[n] = \{1,2,\ldots,n\}$, where $n \in \mathbb N$ is a finite number. A subset $A \subseteq [n]$ is nothing more than a set containing integers between $1$ and $n$, e.g., $A=\{2,4,6\} \subset [7]$.

A family $\F \subseteq 2^{[n]}$ is a collection of subsets of $[n]$. Here $2^{[n]}$ contains all $2^n$ possible subsets of $[n]$, which includes the empty set $\emptyset$ as well.

\begin{figure}[h]
    \centering
    \begin{tikzpicture}[scale=1.5]

\node (c) at (0,3) {$\{1,2\}$};
\node (f) at (0,4) {$\{1,2,3\}$};
\node (g) at (0,5) {$\{1,2,3,4\}$};
\node (a) at (-1,3) {$\{1,3\}$};
\node (e) at (0.5,2) {$\{2\}$};
\node (b) at (-0.5,2) {$\{1\}$};
\node (d) at (1,3) {$\{2,3\}$};

\draw  (g) -- (f) -- (a) -- (b);
\draw (e) -- (c) -- (f) ;
\draw (f) -- (d) -- (e);
\draw (f) -- (c);
\draw (c) -- (e);
\draw (c) -- (b);

\end{tikzpicture}
    \caption{Example of union-closed family}
    \label{fig:UCfam}
\end{figure}
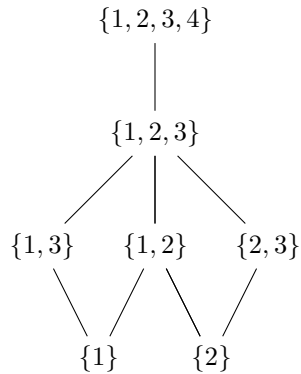

A family $\F$ is called \textbf{union-closed} if for every $A,B \in \F$, the union $A \cup B$ belongs to $\F$. This can be written as $\F=\F \cup \F$, where the latter equals exactly $\{A \cup B \mid A,B \in \F\}.$ An example of such a family is presented in Figure~\ref{fig:UCfam}. An other example, for every $m \in \mathbb N$, is the family $\F_m=\{A \mid A \subseteq [m] \vee A=[k] \mbox{ for some } m+1 \le k \le m^2 \}$ which consists of the $2^m$ subsets of $[m]$, as well as $m^2-m$ intervals consisting of the first $k$ natural numbers.

The Union-closed conjecture can now be formally stated as follows.

\begin{conj}[Union-closed conjecture]\label{conj:UC}
If $\F \neq \{\emptyset\}$ is a union-closed family with ground set $[n]$, then there exists an element $i \in [n]$ such that at least half of the sets in $\F$ contain $i$.
\end{conj}

Considering our previous example $\F_m$ for large $m$, one can verify that it might be that only a small fraction of the elements of the ground set are abundant (belong to at least half of the sets) and their average proportion of sets to which they belong can tend to zero.
Note that this conjecture would be (arguably) false when taking an infinite ground set $\mathbb N,$ e.g. by considering
the (union-closed) family of finite subsets of $\mathbb N$.

This conjecture can also be formulated in many different ways. For example, one can consider bitstrings in $\{0,1\}^n$ with the element-wise $OR$-operation. For instance, when $n=4$ and $\F=\{0011,1100,1111\}$, we note that $0011+1100=1111$. This family is closed under the $OR$-operation, which corresponds to being union-closed in the initial formulation.

Taking the complements of the set, one obtains the Intersection-closed sets conjecture, which states that an intersection-closed family has an element in its ground set appearing in at most half of the sets. In~\cite[Sec.~3]{BS15}, one can also find a lattice-, graph-, and Salzborn-formulation.

On November 17, 2022, Justin Gilmer~\cite{Gilmer22}, a researcher at Google working in machine learning, made a breakthrough by proving a first constant fraction for Conjecture~\ref{conj:UC}. Soon thereafter, as fast as a few days, his result made others put improvements and related results on the preprint server Arxiv.
In this note, we summarize the contributions and progress that was made in the winter of 2022--2023. We explain the main ideas of Gilmer's approach (Section~\ref{sec:proofGilmer}), mention the forthcoming extensions of his method (Sections~\ref{sec:immedfollowup} and~\ref{sec:followup}), as well as an unsuccessful attempt (Section~\ref{sec:eureka}) and discuss other work related to the Union-closed conjecture (Section~\ref{sec:otherdirections}). 

\section{The observations and key elements in the proof by Gilmer}\label{sec:proofGilmer}


A first elementary observation by Gilmer is that one can always prove a statement by proving the contrapositive of that statement. Since the statement of the union-closed conjecture is that simple already, it might be no one considered that before.
The contraposition of Conjecture~\ref{conj:UC} can be stated as follows.
If a non-empty family $\F$ has no element appearing in at least half of the sets of $\F$, then $\F$
is not a union-closed family. 
By remarking that $A \cup A=A$ for every set $A$, one knows that $\F \subseteq \F \cup \F,$ and thus $\abs { \F \cup \F} > \abs \F$ whenever $\F$ is not a union-closed family.
While posing related questions and studying counterexamples to variants of Conjecture~\ref{conj:UC} similar to the ones in~\cite{Ellis20}, Gilmer noted that the entropy of a family might play a role.\footnote{More details on his journey/ thought process can be found in \url{https://www.youtube.com/watch?v=AZaP0EwjR_I&t}}
The entropy $H(X)$ of a discrete random variable $X$ equals the Shannon entropy of its probability distribution. 
The latter can be purely presented with a formula.
If each possible outcome $x$ belongs to a (finite) set $A$, and has probability $p_x,$ then 
$$H(X)=-\sum_{x \in A} p_x \log_2 p_x.$$
When sampling uniformly at random from $\F$, the entropy will equal $\log_2 {\abs \F}$ and no higher entropy is possible.
If one can sample from $\F \cup \F$ in such a way that the entropy is larger than $\log_2 {\abs \F},$ then one can conclude that $\abs { \F \cup \F} > \abs \F.$ 
This is exactly the core of Gilmer's approach.

More precisely, he proved the following statement.

\begin{theorem}\label{thm:Gilmer1}
Let $A$ and $B$ denote independent and identically distributed random variables that sample from a common distribution over subsets of $[n]$. Assume that for all $i \in [n]$, $\Pr[i \in A] \leq 0.01$. Then $H(A \cup B) \geq 1.26 H(A)$.
\end{theorem}

As a corollary, by taking the uniform distribution over the subsets of $[n]$, one knows that if $\F \subset 2^{[n]}$ is a family for which every element is contained in no more than $1\%$ of the sets, then $\abs{ \F \cup \F} \ge \abs{\F}^{1.26}$.\footnote{As a corollary of later work by Sawin, this is at least $\abs{\F}^{1.74}$} This implies that whenever $\abs{\F} \ge 2$, either $\abs{ \F \cup \F} > \abs{\F}$ (and so the family is not union-closed) or there is an element appearing in at least a $0.01$ fraction of the sets in $\F.$ From this, one can conclude that Conjecture~\ref{conj:UC} is true for a half replaced by $0.01.$

\begin{ex}
    Let $\F=\{ \{1\}, \{2\} \} $ and thus
        $\F\cup \F=\{ \{1\}, \{2\}, \{1,2\} \}.$
    Let $A$ and $B$ be i.i.d. random variables that output a set of $\F$ uniformly at random.
    Then $\Pr(A=\{1\})=\Pr(A=\{2\})$ and analogously for $B$, which implies $$\Pr( A \cup B =\{1\}) = \Pr( A \cup B =\{2\}) = \frac 14 \mbox{ and } \Pr( A \cup B =\{1,2\}) = \frac 12.$$
    Now $H(A)= 2 \cdot \frac 12 \log_2 2 =1$ and $H(A \cup B)= 2 \cdot \frac 14 \log_2 4 + \frac 12 \log 2 =\frac 32 (< \log_2 3).$
    Since $\log_2(2)<H(A \cup B),$ we conclude that it is impossible that $A\cup B$ takes values in a family with only $2$ elements and thus $\abs{ \F \cup \F } > \abs \F,$ i.e. Gilmer's method verifies that $\F$ is not union-closed.
\end{ex}

\begin{ex}
    Let $\F=\binom{[3]}{\le 2}$ and thus $\F\cup \F=2^{[3]}.$ Note that $\abs{\F}=7$ and every $1 \le i \le 3$ appears in exactly $3$ sets and thus in a $\frac 37$ fraction.
    Let $A,B$ be i.i.d. random variables that output a set of $\F$ uniformly at random.
    Then
    \begin{align*}
       \Pr(A \cup B = \emptyset) &=\frac 1{49}\\
       \Pr( \abs{A \cup B} = 1) &=\frac 3{49}\\
       \Pr( \abs{A \cup B} = 2) &=\frac 9{49}\\
       \Pr(A \cup B = [3]) &=\frac {12}{49}
   \end{align*}
   Now \begin{align*}
    H(A)=& 7 \frac 17 \log_2(7)=\log_2(7)\\ \sim&2.81\\
        H(A\cup B)=& \frac 1{49} \log_2(49)+3\frac 3{49} \log_2(49/3)  \\ &+3\frac 9{49} \log_2(49/9) +\frac{12}{49} \log_2(49/12) \\ \sim& 2.70
    \end{align*} 
    and thus $H(A)>H(A \cup B)$. We conclude that this is an example for which Gilmer's method does not provide evidence that the family is not union-closed, even while the maximum fraction of occurence of an element is $\frac 37.$
\end{ex}

Note: Analogously, when $\F= \binom{[5]}{ \le 3}$, one can verify that $H(A)=\log_2(26)\sim 4.7$ and $H(A \cup B)\sim 4.54.$
Every element appears in a $\frac{11}{26}$ fraction in this case.

\section{Quick refinement of Gilmer's idea }\label{sec:immedfollowup}

The binary entropy function $h(p)=-( p\log_2 p + (1-p) \log_2(1-p))$ plays a role in the computations in the work of Gilmer. Noting that $h(p) \le h(2p-p^2)$ 
whenever $p \le \psi:=\frac{3- \sqrt 5}{2}$, Gilmer claimed that his ideas could be extended to prove a fraction equal to $\psi.$
The authors of~\cite{AHS22, CL22, Sawin22, Pebody22} quickly implemented this approach. All four of these papers essentially reduced Conjecture~\ref{conj:UC} for the constant $\psi$ to the following key lemma, an inequality in one variable.

\begin{lemma}\label{lem:keylemma_boundpsi}
 Let $\phi=\frac{\sqrt 5+1}{2}$ and $0 \le x \le 1,$ then $h(x^2) \ge \phi x h(x).$
\end{lemma}

The validity of this lemma was established in two different ways by~\cite{AHS22} and Sawin~\cite{Sawin22}. The former used accurate computer calculations and applied interval arithmetic on three intervals, while the latter utilized a purely calculus-based approach. Thanks to some communication between the authors of~\cite{AHS22} and~\cite{CL22}, in~\cite{CL22} a reference to the formal proof of~\cite{AHS22} was added.
In~\cite{Pebody22} the lemma was split in two parts without formal proof, but both can be verified easily.

A short and more elegant proof for Lemma~\ref{lem:keylemma_boundpsi} was given later by Boppana~\cite{Boppana23}, even while the proof itself would originate from $1989$. This proof relies on the following extension of the classical Rolle's theorem, which follows from observations in e.g.~\cite{Petrov03}.

\begin{theorem}
 Let $f$ be a differentiable function on a interval $I$.
 Let $m(f)$ be the sum of multiplicities of the roots of $f$ in $I$.
 Then $m(f') \ge m(f)-1.$
\end{theorem}

By iterating the theorem three times, one finds $m(f)\le m(f''') +3$. Applying this result on the function $f(x)=h(x^2)-\phi x h(x)$ and counting the multiplicities of the roots $0, \frac 1 {\phi}$ and $1$ of $f$, the conclusion that $f$ is nonnegative on $[0,1]$ follows quickly.
Once Lemma~\ref{lem:keylemma_boundpsi} is derived, the proof for Conjecture~\ref{conj:UC} for constant $\psi$ (instead of $0.5$) is rather short in each of the papers~\cite{AHS22,CL22,Pebody22,Sawin22}, indicated e.g. by the total length of the paper by Chase and Lovett~\cite{CL22}. Their work has three steps.
First, they extended the analytic claim (Lemma~\ref{lem:keylemma_boundpsi}) to the two-variate function 
$f(x,y):= \frac{h(xy)}{h(x)y+h(y)x}$.
Next they prove a strengthened inequality between the entropy of $A \cup B$ and the one of $A$ and $B$, for random variables $A$ and $B$ (not necessarily identical) on $\{0,1\}^n$ for which every bit is $1$ with a bounded probability. Finally, they finish the proof of their slightly more general statement that holds for approximate union-closed families. The latter being families for which the union of two random drawn sets belong to the family with a high probability.

One example which certifies the sharpness of their proof can be derived from $\F_1 + \F_2 = \{A \mid A \in \F_1 \vee A \in \F_2\}$ where 
$\F_1=\binom{[n]}{\psi n +n^{2/3}}$ and $\F_2=\binom{[n]}{\ge (1-\psi) n}$. For this, one need to note that $\abs{\F_1} >> \abs{\F_2}$ and that the union of two (iid uniform sampled) random sets from $\F_1$ belongs with very high probability to $\F_2.$ 
The expected size of the union is slightly larger (with an additional term of the order $n^{2/3}$, i.e. $\Theta(n^{2/3})$) than $n-(1-\psi)^2 n=(1-\psi) n$, and since the variance on the size is $O(n^{1/2})$, the union almost surely belongs to $\F_2$ as well.
The conclusion is still valid when replacing the term $n^{2/3}$ by any function $g(n)$ for which $n>>g(n)>>n^{1/2}.$

\begin{figure}[h]
 \centering
 \begin{tikzpicture}[scale=5]
 \def\psi{0.382}
 \def\n{8}
 \def\bigset{{\n-1}}
 \draw[line width=0.05mm, color = lightgray] (-\psi,\psi) -- (\psi,\psi) node at (0,0.45){$\mathcal{F}_1$};
 \draw[line width=0.005mm, fill=lightgray] (-\psi,1-\psi) -- (\psi,1-\psi)--(0,1)--cycle node at (0,0.8) {$\mathcal{F}_2$};
 \draw (-0.5,0.5) -- (0,1) -- (0.5,0.5) -- (0,0) -- cycle;
 \end{tikzpicture}
 \caption{An approximate union-closed family whose elements appear in at most a $\psi +o(1)$ fraction. }
 \label{fig:CLconstr}
\end{figure}
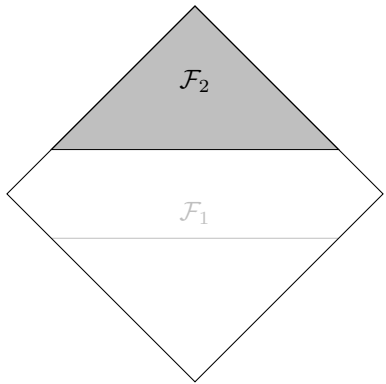

In a different direction, in his paper, Gilmer included some ideas for a full resolution of Conjecture~\ref{conj:UC}, but some of these directions were immediately proven not to hold by Sawin and Ellis~\cite{Sawin22, Ellis22}.

\section{Further refinements and extensions related to Gilmer's work}\label{sec:followup}

Sawin~\cite{Sawin22} gave a suggestion to improve the bound further, which given the sharpness of the form for union-closed families may be considered surprising. Hereby the essence is in a question purely stated in terms of probability distributions.
His suggestion was worked out by Yu~\cite{Yu22} and Cambie~\cite{Cambie22}.
Yu~\cite{Yu22} considered the approach in a slightly more general form initially and made a lower bound computable by restricting to the suggestion of Sawin and applying~\cite[Lem.~5]{AHS22} and the Krein-Milman theorem~\cite{KM40} to bound the support (number of values with nonzero probability) of a joint distribution by $4$.
A numerical computation then yield a bound equal to (roughly) $0.38234.$
In parallel, Cambie~\cite{Cambie22} found an upper bound for Sawin's approach which indicates that the improvement is way smaller than expected and one would hope for. The construction is a discrete probability distribution with only two values having nonzero probability, with the values determined by a system of equations involving the entropy function.
Additionally he proved that this value is sharp, by first reducing the support to $3$ elements, where one of the elements equals $1.$
Finally, the conclusion is derived from the combination of $3$-dimensional plots, a numerical minimization problem and a more precise solution for the case where the support has exactly two elements, one of which equals $1$.

Finally, building upon the work of~\cite{CL22}, Yuster~\cite{Yuster23} considered families that are almost $k$-union-closed, meaning that the union of $k$ independent uniform random sets from $\F$ belongs to $\F$ with high probability. He conjectured a tight version for the minimum frequency (the proportion of sets containing the element) of some element in such families, with the threshold for this frequency being the unique real root in $[0,1]$ of $(1-x)^k=x$, denoted by $\psi_k$. To understand the sharpness of his conjecture and the intuition behind the choice of $\psi_k$, consider the union of $\F_1=\binom{[n]}{\psi_k n^+n^{2/3} }$ and $ \F_2=\binom{[n]}{\ge (1-\psi_k) n}$. If at least one set from $\F_2$ is included among the $k$ sets drawn, the union is guaranteed to belong to $\F_2$. If all $k$ sets belong to $\F_1$, the expected size of the union is $n-(1-\psi_k)^k n+\Theta(n^{2/3})$, and since the variance is $O(n^{1/2})$, the union almost surely belongs to $\F_2$ as well. The conjecture is proven to be true for $k\leq 4$, while for larger values of $k$ a weaker bound is established.

\section{The final Eureka moment, not yet}\label{sec:eureka}

When Scandone~\cite{Scandone23} uploaded a preprint claiming the full resolution of the union-closed conjecture, there arose initially excitement. 
However, upon closer examination it became clear that Scandone's proposed solution had several issues, including a significant flaw that requires revising the underlying construction.
This was communicated to Scandone by Terence Tao, and the details of this issue are briefly explained later in this section.

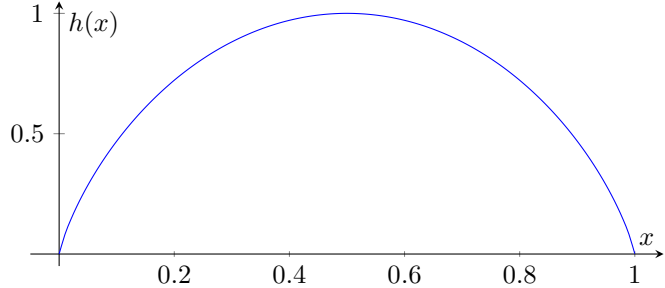
\begin{figure}[h]
 \centering
 \begin{tikzpicture}
\begin{axis}[
domain=0:1,
samples=100,
yscale=0.5,
xlabel={$x$},
ylabel={$h(x)$},
axis lines=middle,
xtick={0,0.2,...,1},
ytick={0,0.5,...,1},
xmin=-0.05, xmax=1.05,
ymin=-0.05, ymax=1.05,
clip=false
]
\addplot[blue, smooth] {1/ln(2)*(-x*ln(x) - (1-x)*ln(1-x))};
\end{axis}
\end{tikzpicture}
 \caption{Plot of the binary entropy function $h$}
 \label{fig:plot_h}
\end{figure}

Nevertheless, Scandone's underlying idea holds potential and is worth mentioning for the valuable intuition it provides for Gilmer's approach.
Let $\F$ be a family which is not union-closed, so $\F \cup \F \not= \F$.
A random variable taking values in $\F$ has entropy at most $\log_2 \abs \F$ and equality occurs only for uniform sampling from $\F.$
By considering various examples, e.g. $\F=\{ \{1\},\{2\} \}$, the reader can verify that 
there is no strategy to choose two random variables $A,B$ which sample sets from $\F$, such that $A\cup B$ samples uniformly random from $\F \cup \F$.
On the other hand, if for every set $A \in \F$ the probability of obtaining it is almost equal to the original probability and a few other sets from $(\F \cup \F) \backslash \F$ happen with a small probability, the entropy can increase.
The reason for this is that the derivative of $h$ (plotted in Figure~\ref{fig:plot_h}) is a continuously decreasing function on the interval $(0,1)$, with $h'(0)= + \infty$.
To provide a more explicit explanation of Scandone's idea, we describe his proposed construction in detail.

%

Let $A, B$ be independent random variables that take any set of $\F$ uniformly at random.
Define a $\mathcal{P}([n])$-valued random variable $A^{\delta}$ (depending on $\delta$) through the relation
$$\pr[A^{\delta}=X]=(1-\delta)\pr[A=X]+\delta\pr[A\cup B=X] \mbox{ for every } X \subseteq [n].$$
For every $X \in \F$, $\pr[A^{\delta}=X]\ge (1-\delta)\pr[A=X]$ and thus for $\delta$ sufficiently small, we have
$h(\pr[A^{\delta}=X])-h(\pr[A=X]) \gtrsim \delta / \abs \F h'( 1/ \abs \F).$\footnote{To be precise, we assume $\abs \F \ge 3$ and $\frac 2{\abs \F} + \delta <1.$}
On the other hand, for $X \in (\F \cup \F) \backslash \F$, let the probability $p:= \pr[A \cup B=X]$.
We have that $h( \delta p) \sim -\delta p (\log \delta + \log p-1)$.
By choosing $\delta$ to be sufficiently small such that $-\log \delta$ is much greater than $\frac{1}{p}h'(1/\abs{\F})$, we can ensure that $H(A^\delta)>H(A)$ holds.

Equivalently, the variable $A^{\delta}$ can be obtained by considering, in addition to $A$
and $B$, a Bernoulli random variable of parameter $\delta$, $Z_{\delta}$, which determines whether we take $A \cup B$ or only $A$. The flaw in the argument is that, in
the process of revealing all the digits of $A^{\delta}$ (computed using the chain rule for the entropy), the indeterminacy
provided by $Z_{\delta}$ (and the consequent improvement of the bounds) is lost after the first step. 
More precisely, there is step in the computations in which
a conditional probability distribution has been erroneously replaced by
its expected value, and this produces the aforementioned flaw in the
argument.
The comment of Tao can be rephrased as follows, ``the idea
of modifying the union operation by Gilmer is promising, but a single
global bit $Z_{\delta}$ is not sufficient to do the job, and a more
involved construction is needed''.

\section{A better understanding by progress in a different direction}\label{sec:otherdirections}

In this final section, we conclude with the essence of a recent paper and two preprints on the union-closed conjecture, which consider different aspects and angles of attack on Conjecture~\ref{conj:UC}.

While Frankl's conjecture is about the existence of one abundant element (element that appears in at least half of the sets) in the family, it is also natural to wonder if there are more abundant elements, assuming that all sets in the family are sufficiently large.
The following conjecture by Cui and Hu~\cite{CuiHu21} would imply Conjecture~\ref{conj:UC}.
\begin{conj}\label{conj:CuiHu}
 If $\F$ is a finite union-closed family of sets whose smallest set is of size at least $2$,
then there are at least two elements such that each belong to more than half of the sets of $\F$.
\end{conj}
At the end of $2022$, the three authors of~\cite{KPT22} considered this different direction and proved that Conjecture~\cite{CuiHu21} is not true when replacing $2$ by a larger integer.
They proved (among other results) that there are families all of whose sets have size at least $k$, where $k$ can be arbitrary large, which do only have $2$ abundant elements.
The main construction is the family $\Ss^{12}_4$.
The family $\Ss^{12}_4$ consists of all subsets $S$ of $\{0,1,\dots,11\}$ of size at least $4$ such that 
either $\{0,1\} \subset S$, or
$0 \in S$ and $S \subseteq \{0,2,\dots,10\}$, or
$1 \in S$ and $S \subseteq \{1,3,\dots,11\}$.
The reader can verify that $\abs{ \Ss^{12}_4} = (2^{10}-11)+2\cdot 16=1045$, while every element $2 \le i \le 11$ only appears $2^9-1 +11=522$ times.
One way to increase the size of sets in families with non-abundant elements is to duplicate an element within the sets. However, this creates blocks of size at least $2$. A block is defined by Poonen~\cite{Poonen92} as a maximum set of elements that all belong to the exact same sets of a family. Poonen also noted that to prove Conjecture~\ref{conj:UC}, it is sufficient to focus on families for which no block is a singleton.
Due to this, it is interesting to note that the construction of the family $\Ss^{12}_4$ in~\cite{KPT22} can be extended to such families.
Let $k\ge 3$ be a fixed integer and let $n$ be a sufficiently large even integer as a function of $k$ ($n\ge 10k$ works).
Let $E_n=\{i \in [n] \mid i \equiv 0 \pmod 2\}$ and $O_n=\{i \in [n] \mid i \equiv 1 \pmod 2\}$ be the set of even and odd integers in $[n]$ respectively.
Consider the family $\Ss^{n}_k$ consisting of subsets $S$ of $[n]$ of size at least $k$, such that either
\begin{itemize}
 \item $\{1,2\} \subset S$,
 \item $S \subset E_n$ and $2 \in S$, or
 \item $S \subset O_n$ and $1 \in S$.
\end{itemize}
It is clear that $1$ and $2$ are abundant elements.
Now the other elements appear all equally often (by symmetry) and by a small bijection and counting argument, we conclude that these elements are not abundant whenever
$$ \binom{n-3}{k-3} < 2 \binom{n/2-2}{\ge k-1}.$$
Since this is the case for $n$ sufficiently large, the conclusion is clear.

Another result related with union-closed families and the smallest set size, was published early $2023$.
Ellis, Ivan and Leader~\cite{ELM23} proved that for every $k \in \mathbb N$, there exists a union-closed family in which the (unique) smallest set has size $k$, but where each element of this set has frequency
$(1+o(1))\frac{\log k}{2k}.$ As such, proving that focusing on the smallest set cannot work in the strongest possible sense. They also proposed the problem of verifying the union-closed conjecture for a family for which they were unable to verify the statement.
The latter was verified by Pulaj and Wood~\cite{PW23}. They also proved new bounds on the least number $m$ (given $k$ and $n$) such that every union-closed family $\F$ containing any $\mathcal A \subseteq \binom{[n]}{k}$ with $\lvert \mathcal A \rvert =m $ as a subfamily, satisfies Conjecture~\ref{conj:UC}.

We can conclude that despite the progress that originates from the breakthrough of Justin Gilmer, the exact version of Conjecture~\ref{conj:UC} is still not proven.
Mathematicians are still thinking about other directions or modifications of the strategy and hope to resolve Conjecture~\ref{conj:UC} in the future.
Taking into account that the improvement by taking combinations suggested by Sawin~\cite{Sawin22} turned out to be tinier than expected and hoped for, as illustrated by the example in~\cite{Cambie22}, it seems that the focus should go towards essential new ideas.
In particular, the union-closed conjecture might be a distraction of a more general behaviour that $\abs{\F \cup \F} > \abs{\F}^c $ for some $c(\eps)>1$ when every element of $[n]$ appears in less than a $\frac 12- \eps$ fraction of the sets in $\F.$\footnote{communicated by Zachary Chase} 

Note added: In June 2023, Liu~\cite{Jingbo23} improved the constant slightly with a different method of coupling.

\section*{Acknowledgements}
We thank Zachary Chase, Justin Gilmer, Raffaele Scandone and Lei Yu for internal communication while writing this manuscript.

\bibliographystyle{abbrv}
\bibliography{unionclosed}

\end{document}